  \documentclass{amsart}

\usepackage{amsmath,amssymb,amscd,amsfonts,verbatim}

\usepackage[mathscr]{eucal}

\newtheorem{thm}{Theorem}[section]

\newtheorem{lem}[thm]{Lemma}
\newtheorem{rem}[thm]{Remark}
\newtheorem{prop}[thm]{Proposition}
\newtheorem{cor}[thm]{Corollary}

\newtheorem{assu-nota}[thm]{Assumption--Notation}
\theoremstyle{remark}

\newcommand{\pionealg}{\pi_1^{\mathrm{alg}}}

\newcommand{\C}{\mathbb C}
\newcommand{\Z}{\mathbb Z}

\newcommand{\pp}{\mathbb P}

\newcommand{\Ical}{\mathcal{I}}

\newcommand{\Rcal}{\mathcal{R}}
\DeclareMathOperator{\Aut}{Aut}
\DeclareMathOperator{\Pic}{Pic}

\DeclareMathOperator{\Alb}{Alb}

\newcommand{\epsi}{\varepsilon}

\newcommand{\ga}{\gamma}

\newcommand{\Ga}{\Gamma}
\newcommand{\De}{\Delta}
\newcommand{\Si}{\Sigma}

\newcommand{\fie}{\varphi}

\newcommand{\OO}{\mathcal{O}}
\newcommand{\Ocal}{\mathcal{O}}
\newcommand{\inv}{^{-1}}

\numberwithin{equation}{section}

\title{Surfaces with $K^2<3\chi$ and finite fundamental group}
\author{Ciro Ciliberto,  Margarida Mendes Lopes and Rita Pardini}

  \address{Ciro Ciliberto\\ Dipartimento di Matematica\\ Universit\`a
degli Studi di Roma ``Tor Vergata''\\Via della Ricerca Scientifica, I-00133 Roma, Italy }

  \email{cilibert@mat.uniroma2.it}

\address{Margarida Mendes Lopes\\Center for Mathematical Analysis, Geometry and Dynamical Systems\\
Departamento de  Matem\'atica\\ Instituto Superior T\'ecni\-co,  
Universidade T\'ecnica de Lisboa\\
Av.~Rovisco Pais, 1049-001 Lisboa, Portugal}

\email{mmlopes@math.ist.utl.pt}

\address{Rita Pardini\\Dipartimento di Matematica\\
Universit\`a di Pisa\\
Largo B. Pontecorvo 5,
56127 Pisa, Italy}

\email{pardini@dm.unipi.it}

\begin{document}
\begin{abstract}

 In this paper we continue the study of  $\pionealg(S)$ for minimal surfaces of general type  $S$ satisfying $K_S^2 <3\chi(S)$. We show that, if $K_S^2= 3\chi(S)-1$ and $|\pionealg(S)|= 8$, then
  $S$ is a Campedelli surface. 
  
  In view of the results of \cite{3chi} and \cite{3chi-2}, this implies  that the fundamental group of a surface   with $K^2<3\chi$ that has no irregular \'etale cover has order at most 9,  and if it has order 8 or 9, then $S$ is a Campedelli surface.

To obtain this result we establish some classification results for minimal surfaces of general type such that $K^2=3p_g-5$ and such that the canonical map is a birational morphism.  We also study rational surfaces with a $\Z_2^3$-action.
\end{abstract}
\maketitle

\section{Introduction}

The algebraic fundamental group of a surface of general type is an important invariant of the surface. If a surface is irregular then of course $\pionealg$ is infinite, but regular surfaces may also have infinite fundamental group.
 
For minimal surfaces $S$ of general type with $K^2< 3\chi$ through the work of several authors (\cite{bombieri}, \cite{ho5}, \cite{miles1}, \cite{milesk2}, \cite{xiaohyp}, \cite{xiaoslope}, \cite{mur1}, \cite{mur2}, \cite{cm}, \cite{3chi}, \cite{3chi-2}) one has  a quite good picture of  what can occur, and in particular a very good understanding of the geometry of $S$ when $S$ has infinite fundamental group (cf. the introduction of \cite{3chi}). It turns out that $S$ has an infinite fundamental group if and only if it has an irregular \'etale cover. Also, if $S$ has finite fundamental group the groups that can occur form a very limited family. In particular if $K^2\leq 3\chi- 2$, then the order of $\pionealg$ is at most $5$ and if $|\pionealg(S)|=5$ then $S$ is a Godeaux surface (i.~e. $K^2=1,  p_g(S)=0, \chi(S)=1$) (see \cite{3chi-2}). If $K^2=3\chi-1$, then the order of $\pionealg$ is at most $9$ and if $|\pionealg(S)|=9$ then $S$ is a Campedelli surface (i.~e. $K^2=2, p_g(S)=0, \chi(S)=1$).

In this paper we sharpen this last result and we prove   the following:

\begin{thm} \label {thm:m} Let $S$ be a minimal complex surface of general
type such that $K^2_S=3\chi(S)-1$. If the group $\pionealg(S)$ has order
8,
then $\chi(S)=1$.  \end{thm}

By the results in \cite{3chi, 3chi-2}, we have then the following:

\begin{thm}\label{main2} Let $S$ be a minimal complex surface of general
type with $K_S^2\leq 3\chi(S)-1$ which has no irregular \'etale cover. Then
the
group $\pionealg(S)$ is finite of order $\le 9$.

If $|\pionealg(S)|=8$ or 9, then $\chi(S)=1$, $K^2_S=2$, namely $S$ is a
numerical Campedelli surface.

 \end{thm}

\medskip Note that there exist Campedelli surfaces with algebraic fundamental group of orders 8 and 9 (cf. \cite{milesk2} and \cite{camp9}). 

\medskip

The idea for the proof of Theorem \ref {thm:m} is the following. Let $V\to
S$ be the \'etale cover of degree 8 given by $\pionealg(S)$. Then $V$ is
regular and $K_V^ 2=3p_g(V)-5$, i.e. these invariants are {\it close} to
the {\it Castelnuovo line} $K^ 2=3p_g-7$.

Canonical surfaces with $3p_g-7\leq K^ 2\leq 3p_g-6$ have been studied and
classified in \cite{harris}, \cite{ak} and \cite{konnopg}. The case $K^
2=3p_g-5$, of interest here, is however still widely open. In \S \ref
{sec:cansurf} we start the analysis of canonical surfaces with these
invariants and obtain the information we need for  the proof of Theorem \ref
{thm:m}. Namely, we consider a surface $V$ with $K_V^ 2=3p_g(V)-5$, for
which the canonical map is a birational morphism of $V$ onto its image
$\Sigma$. One proves that $V$ is regular (see Lemma \ref{q}) and that, if
$p_g\geq 14$, then $\Sigma$ lies on a rational threefold scroll $X$ of
minimal degree, whose ruling pulls back to a base point free pencil of
curves of genus 3 on $V$ (see Theorem \ref {g3} for a more precise
statement). This result opens up the possibility for a
detailed classification of these canonical surfaces. However we do not
dwell on this here, since Theorem  \ref {g3} is what we need   for our  proof of Theorem \ref {thm:m}. The proof of the results of \S \ref
{sec:cansurf} use Castelnuovo's theory, in particular the geometric
information  on a set of points in uniform position which impose {\it few}
conditions on quadrics (see \cite {h}, \cite {eh}, \cite {c}).

Next, having in mind the application to the proof of Theorem \ref {thm:m},
we consider, in \S \ref {nonbir}, the following general situation: $S$ is
a
minimal surface of general type, $V\to S$ is an \'etale  Galois cover with
Galois group $G$ of order 8, the canonical map of $V$ has degree 2 onto a
rational surface $W$. The canonical involution lies in the centre
of ${\rm Aut}(V)$, hence the action of $G$ on $V$ descends to an action of
$G$ on $W$.  Moreover, a result by Beauville  \cite {beauville}
implies that $G=\Z_2^3$. We analyse this situation by considering the
{\it relative minimal model} \`a la Mori of the pair $(W,G)$. We show  that it is a conic bundle and that the subgroup $H$ of $G$ that sends each fibre to itself has order $2$. Finally we describe in detail  the whole situation (see Proposition \ref {minimal}).

Finally, in \S \ref {sec:fund} we go back to the problem we started with.

 If $S$ has $K^ 2_S=3\chi(S)-1$ with $G=\pionealg(S)$ of order
8, and $V\to S$ is the \'etale cover given by $G$, we prove that either
the
canonical map of $V$ is a birational morphism to its image, or it has
degree 2 onto a rational surface (see Proposition \ref {canV}). Then, by
using the aforementioned results of \S\S   \ref{sec:cansurf}, \ref {nonbir}, we separately discuss the two cases, excluding
the
latter, and  proving that the former is compatible only with
$\chi(S)=1$, i.e. $S$ has to be a Campedelli surface.

\bigskip

{\bf Notation and conventions.} We work over the complex numbers. All
varieties are projective algebraic. All the notation we use is standard in algebraic geometry. 
We just recall the definition of the nume\-rical invariants of a smooth surface $S$:  the self-intersection number $K^2_S$ of the canonical divisor $K_S$, the {\em geometric genus} $p_g(S):=h^0(K_S)=h^2(\OO_S)$, the {\em irregularity} $q(S):=h^0(\Omega^1_S)=h^1(\OO_S)$ and the {\em holomorphic Euler characteristic} $\chi(S):=1+p_g(S)-q(S)$.

\section{Preliminaries}
We recall some general facts that will be used in the following sections.
\subsection{Fixed point formulae}\label{ssTF}
Let $S$ be a smooth complex compact surface and let  $\iota$ be an involution of $S$, namely an automorphism of order 2 of $S$.  The fixed locus of $\iota$ is the union of a smooth (possibly empty) divisor $R$ and of  $\nu$
isolated points. We
recall
 the following
two well-known formulae:

Holomorphic Fixed Point Formula (see \cite{as}, pg. 566):
\begin{equation}\label{HTF}
\sum_{i=0}^2(-1)^i\text{Trace}(\iota|H^i(S,\OO_S)) = \frac{\nu-RK_S}{4}
\end{equation}
Topological Fixed Point Formula (see \cite{gre}, (30.9)):
\begin{equation}\label{TTF}
\sum_{i=0}^4(-1)^i\text{Trace}(\iota|H^i(S,\C)) = \nu+e(R),
\end{equation}
where $e(R) = -R^2-RK_S$ is the topological Euler
characteristic of $R$.

If $S$ is a rational surface,  as it is often the case in this paper,  then (\ref{HTF}) takes a particularly simple form:
\begin{equation}\label{TFrat}
\nu=4+K_SR
\end{equation}
If, in addition, the divisor $R$ is empty, then (\ref{TFrat}) gives $\nu=4$.
\begin{rem}\label{remsing}{\rm Formulae \eqref{HTF} and \eqref{TFrat}  are stated here  for a smooth surface, but clearly
they hold also for a  surface with rational singularities,   provided that the involution acts
freely on the singular locus.}
\end{rem}
\subsection{Free group actions on pencils of low genus} We will use the following result from \cite{3chi-2}.
\begin{lem}\label{free} Let $Y$ be a regular surface of general type,  let $G\ne \{1\} $ be  a finite group   that acts freely on $Y$ and  let $|F|$ be a $G$-invariant free pencil $|F|$ of curves of genus $g(F)\leq 4$. Then only the following possibilities can occur:
\begin{enumerate}
\item $G=\Z_2^2$, $g(F)=3$ and $G$ acts faithfully on $|F|$;
\item $G=\Z_3$, $g(F)=4$;
 \item $G=\Z_2$, $g(F)=3$. 
\end{enumerate}
\end{lem}

\subsection{Threefolds  with special curve sections}\label {sec:deltagenus}

 The results of this section are related to  the classification of varieties with small {\it Delta--genus} (cf. \cite{cr}, \cite{eh}, \cite {Fuj1},  \cite {Fuj2}, \cite {Fuj3}, \cite {Fuj4}, \cite {Fuj5}).
 
 Recall that a variety  in a projective space   is a {\it scroll} over a smooth curve $\Delta$ if it is swept out by a 1--dimensional family $\Rcal$ of $n\!-\!1$-dimensional linear subspaces  parametrized by a curve birational to $\Delta$. The elements of the family $\Rcal$ are called {\it rulings}. A linearly normal variety which is a scroll over $\pp^1$ is called a {\it rational normal scroll}. These have been classically described by del Pezzo (see \cite{eh}).

Throughout this subsection we denote by   $X\subset \pp^{r+1}$   a non--degenerate,   irreducible threefold,  by $Y$ its general hyperplane section and by  $\ga$ the arithmetic genus of a general curve section of $X$.  In addition,  we  assume  that $X$ is {\em linearly normal}, namely that  $h^ 1(\pp^{r+1}, \Ical_{X,\pp^ {r+1}}(1))=0$.

We need the following:
\begin{lem}\label {scroll} If  $\gamma>0$  and   $Y$  is a scroll over a curve $\Delta$, then also $X$ is a scroll over $\Delta$.
\end{lem}

\begin{proof} Let $f\colon W\to X$ be a desingularization of $X$ and set  $f^ *(Y)=S$. Note that $S$ is smooth. The divisor $S$ is also nef and big, hence  for $i<3$ one has $H^i(W,\Ocal_W(-S))=0$.  This implies that $h^1(W,\Ocal_W)=h^ 1(S,\Ocal_S)$  and the map $\Alb(S)\to\Alb(W)$ is an isogeny. 

Let $R$ be the proper transform on $S$ of a general line of the ruling of
$Y$. The Albanese map $a_S\colon S\to\Alb(S)$ maps onto $\Delta$ with
general fibre $R$.  Since the
system $|S|$ on $W$ is birational and $S$ is  general, it follows  that
the
Albanese image of $W$ is also a curve $\Delta'$.

Let $\Pi$ be the general fibre of the Albanese map of $W$, which is irreducible. Then the
intersection of $S$ with $\Pi$ is also irreducible and, by the above
argument, it coincides with $R$. This yields that $f$ maps $\Pi$ to a
linear subspace of $\pp^ r$ of dimension $n-1$, proving the
assertion.\end{proof}

 If $\ga\le 2$, then  the possibilities for $X$ are quite restricted:

\begin{prop}\label{threefold}
\begin{enumerate}
\item If  \ $\ga=0$, then $X$ is either a rational normal scroll or $X\subset\pp^6$ is the cone over the Veronese surface in $\pp^5$;
\item if $\ga=1$ and $X$ is not a scroll then   $r\leq 9$;
\item if $\ga=2$  and $X$ is not a scroll  then   $r\leq 11$.
\end{enumerate}
\end{prop}
\begin{proof} Assume $\ga=0$. Then the general curve section of $X$ is smooth irreducible, hence a rational normal curve,  and it is  a classical result (see \cite[Theorem 1]{eh}) that $X$ is as in statement (i).

Assume $\ga>0$ and $X$ is not a scroll. By Lemma \ref{scroll} $Y$ is not a scroll, either. Let $Y_0\to Y$ be the minimal desingularization of $Y$ and let $D$ be the pull-back on $Y_0$ of the general hyperplane section of $Y$. The divisor  $D$ is smooth  and nef, hence statements (ii) and (iii) follow immediately by applying  \cite[Theorem 7.3]{cr} to the pair $(Y_0,D)$.
\end{proof}

\subsection{Castelnuovo theory }
We recall  some  facts from Castelnuovo  theory. For any variety $T\subset \pp^m$  we denote by  $h_T$
the Hilbert function of $T$, namely $h_T(n)$ is the dimension of the image of the restriction map $H^0(\OO_{\pp^m}(n))\to H^0(\OO_T(n))$. 

Given a variety $X\subset \pp^m$, we denote its general hyperplane section by $Y$. Then

\begin{equation}\label {eq:fact1} h_X(n)-h_X(n-1)\geq h_Y(n)
\end{equation}
for all $n\geq 1$. Equality holds in \eqref{eq:fact1} for all $n\geq 1$ if and only if $h^ 1(\pp^m,\
\Ical_{X,\pp^ m}(n))=0$ for all $n\geq 0$, i.e. if and only if $X$ is {\it
projectively normal} (see {\cite[Lemma 3.1 and Remark 3.1.1]{h}; see \cite[ Remark 1.8, (ii)]{c},
 for a more precise result).
\medskip

Of special interest to us here is the Hilbert function of  a finite set $\Ga\subset\pp^{r-1}$. If the points of $\Ga$ are in general position and $d\ge 2r-1$, then it is easy to check that $h_{\Ga}(2)\ge 2r-1$. The cases in which $h_{\Ga}(2)$ is close to this lower bound have been studied by several authors:
\begin{prop}\label{Cast} {\rm  (Castelnuovo, \cite[Lemma 3.9] {h})} Let \,$r\geq 3$ and let $\Gamma\subset \pp^{r-1}$ be a finite  set of $d\geq 2r+1$ points in general position. 
 Then $h_\Gamma(2)=2r-1$ if and only if $\Ga$ is contained in   a rational normal curve of degree $r-1$ in $ \pp^{r-1}$, cut out by all the quadrics containing $\Ga$.
\end{prop}
In order to state the remaining results, we need to recall the following definition (cf. \cite[p.85-86]{h}).
A  finite set $\Ga\subset \pp^{r-1}$ of $d$ points is said to be in {\em uniform position} if, for any subset  $\Ga'\subset \Ga$ of $d'$ points and for any $n\ge 0$,  one has $h_{\Ga'}(n)=min\{d', h_{\Ga}(n)\}$. This definition is motivated by the fact that the general hyperplane section of an irreducible projective curve consists of points in uniform position (cf. \cite[Lemma 3.4]{h}).

\begin{prop}\label{Fano} {\rm (\cite[p.54] {f}; \cite[p.106] {h})} Let \,$r\geq 4$ and let $\Gamma\subset \pp^{r-1}$ be a finite  set of $d\geq 2r+3$ points in uniform position. If  $h_\Gamma(2)=2r$, then  $\Gamma$ is contained  in an elliptic normal curve of degree $r$ in $ \pp^{r-1}$, cut out by all the quadrics containing $\Gamma$.
\end{prop}

\begin{prop}\label{Ciro} {\rm (\cite [Theorem 3.8] {c}, cf. also \cite[Proposition 4.3]{petrakiev})}  Let \,$r\geq 6$, $d>\frac{8}{3}r$ and let $\Gamma\subset \pp^{r-1}$ be a finite  set of $d\geq 2r+1$ points in uniform position. If $h_\Gamma(2)=2r+1$, then  $\Gamma$ is contained  in  an irreducible  curve of degree $\leq r+1$ in $ \pp^{r-1}$, cut out by all the quadrics containing $\Gamma$.
\end{prop}

\section{Canonical surfaces with $K^ 2=3p_g-5$.}\label{sec:cansurf}

Throughout this section we assume  that $V$ 
is  a smooth minimal projective surface of general type  such that:
\begin{itemize}
\item   $p_g:=p_g(V)\geq 4$,  $K^2:=K^2_V=3p_g(V)-5$;
\item  the canonical system  of $|K_V|$ is base point free;
\item the canonical map of $V$ is  birational.
\end{itemize}
 Our analysis is   similar to those of \cite{harris}, \cite{ak} and \cite{konnopg} for canonical surfaces with $K^2=3p_g-7$ and $K^2=3p_g-6$.

Our first remark is that $V$ must be regular: 
\begin{lem}\label{q} One has $q(V)=0$. 
\end{lem}
\begin{proof} 
By \cite[Theorem~3.2]{debarre}, a minimal surface of general type with birational canonical map satisfies: $$K^2\geq 3p_g+q(V)-7.$$
 Therefore, for $V$ satisfying the hypothesis of the lemma, one has $q(V)\leq 2$.
It follows in particular  that $\chi(V)\geq p_g-1$, and therefore $K^2=3p_g-5\le 3\chi(V)-2<3\chi(V)$.
 
 Suppose that $V$ is irregular.    By \cite{ho5}, the condition $K^2<3\chi$ implies that the Albanese   image of $V$  is a curve and that the general fibre $F$ of the Albanese map is an hyperelliptic curve. Since the restriction of  $|K_V|$ to $F$ is contained in $|K_F|$,  the canonical map of $V$ is not birational, contradicting the assumptions.
\end{proof}
\bigskip

 Set $r=p_g-2$, $d=K^ 2=3r+1$ and $g=K^2+1=3r+2$.
\medskip

Let $\Sigma$ denote the canonical image of $V$. Our assumptions imply that  $\Sigma$  is  a non degenerate surface of degree $d$ in $\pp^ {r+1}$.  We denote by $\Delta$ the 
general hyperplane section of $\Sigma$, which is a non--degenerate curve
of degree $d$ in $\pp^ r$. The curve $\Delta$ is linearly normal, since  $q(V)=0$ by Lemma \ref{q}. The curve  $C=\phi^ *(\Delta)$
 is a general curve in $|K_V|$, hence it is  smooth of genus $g$.

We let $\Gamma$ be the general hyperplane section of $\Delta$.  By \cite[Lemma  3.4]{h}, the set  $\Gamma$
consists of $d$ distinct points in uniform position spanning a $\pp^ {r-1}$. As in \S  \ref{sec:deltagenus}, we denote  by $h_\Gamma$ [resp. by $h_\Delta$, $h_\Sigma$] 
the Hilbert function of $\Gamma$ [resp. of $\Delta$, $\Sigma$]. 
From formula \eqref{eq:fact1} one has:
\begin{equation} \label{fact3}  h_\Delta(n)-h_\Delta(n-1)\geq h_\Gamma(n)
\end{equation}
\begin{equation}\label{fact4}  h_\Sigma(n)-h_\Sigma(n-1)\geq h_\Delta(n)
\end{equation}
for all $n\geq 1$.

For any
non--negative integers $n,m$, one also has:
\begin{equation}\label {fact2} h_\Gamma(n+m)\geq \min\{3r+1,h_\Gamma(n)+h_\Gamma(m)-1\}
\end{equation}
(see \cite[Corollary 3.5] {h}).

\begin{lem}\label {lem:hilb} One has:
\begin{enumerate}
\item $h_\Delta(2) \leq g=3r+2$ and, if equality holds, then 
$h^ 1(\pp^ r,\Ical_{\Delta,\pp^ r}(2))=0$ and $\Delta$ is smooth;
\item $2r-1\leq h_\Gamma(2)\leq 2r+1$.
\end{enumerate}
\end{lem} 

\begin{proof} Note that $h_\Delta(2)$ is the dimension of the image of the restriction map

$$H^ 0(\pp^ r,\Ocal_{\pp^ r}(2))\to H^ 0(\Delta,\Ocal_{\Delta}(2))$$
By pulling back to $V$ via the canonical map, this is bounded above by the dimension of the image of the map

$$H^ 0(V,\Ocal_{V}(2K_V))\to H^ 0(C,\Ocal_{C}(2K_V))= H^ 0(C,\omega_C).$$
Since $h^0(C,\omega_C)=g$, assertion (i) follows immediately.
Note now that $h_\Gamma(1)=r$. Hence (ii) follows by (i) and by formulae \eqref {fact3} and \eqref {fact2}. \end{proof}

We first examine the cases $2r-1\leq h_\Gamma(2)\leq 2r$.

\begin {prop}\label {prop:vol1} Assume that $r\ge 4$ (i.e. $p_g\geq 6$).  If $h_\Gamma(2)=2r-1+i$, where  $i=0$ or $1$, then the intersection of all quadrics containing $\Si$ is a linearly normal threefold $X$ of degree $r-1+i$ with curve sections of geometric genus $i$.
\end{prop}
\begin{proof} Assume first $i=0$. Then   Proposition \ref{Cast} (Castelnuovo's Lemma) says that the intersection of all quadrics containing $\Gamma$ is a rational normal curve $Z$ of degree $r-1$.

 Since $q(V)=0$ by Lemma \ref{q}, $\Delta $ is  linearly normal and so the restriction map
$$H^ 0(\pp^ r,\Ical_{\Delta,\pp^ r}(2))\to H^ 0(\pp^ {r-1},\Ical_{\Gamma,\pp^{ r-1}}(2))$$
is surjective.  As a consequence, we have that  the intersection of all quadrics containing $\Delta$ is a surface $Y$ of degree $r-1$, whose  general hyperplane section  is $Z$.

Since $\Sigma$ is linearly normal,
 also the map
$$H^ 0(\pp^ {r+1},\Ical_{\Sigma,\pp^ {r+1}}(2))\to H^ 0(\pp^ {r},\Ical_{\Delta,\pp^ r}(2))$$
is surjective, and, as in the preceding paragraph, we conclude that the intersection of all quadrics containing $\Sigma$ is  a threefold $X\subset \pp^{r+1}$ of degree $r-1$.  The threefold $X$ is linearly normal because $\Sigma$ is also linearly normal.

If $i=1$, then, by Proposition \ref{Fano}, the intersection of all quadrics containing $\Gamma$ is a linearly normal curve $E$ of degree $r$ and genus 1 in $\pp^ {r-1}$.  Then the same argument as above proves the assertion.
\end{proof}

Next we consider the cases $h_\Gamma(2)=2r+1$ and we prove a similar result. 

\begin {prop}\label {prop:vol2}  Assume that $r\ge 6$ (i.~e. $p_g\geq 8$).  If $h_\Gamma(2)=2r+1$,  then the intersection of all quadrics containing $\Si$ is a linearly normal threefold $X$ of degree $r+1$ with curve sections of arithmetic genus either $1$ or $2$.
Moreover,  $\Sigma$ is projectively normal and smooth in codimension 1.
\end{prop}
\begin{proof}  If $h_\Gamma(2)=2r+1$, then we have $h_\Delta(2)=3r+2$ by \eqref {fact3} and by Lemma \ref{lem:hilb}. So, again by Lemma \ref{lem:hilb}, $\Delta$ is smooth and $h^ 1(\pp^ r,\Ical_{\Delta,\pp^ r}(2))=0$. In particular $\Sigma$ is smooth in codimension 1.

  Moreover, one has 
$h^ 1(\pp^ r,\Ical_{\Delta,\pp^ r}(n))=0$, for all $n\geq 3$ by a well known result of 
Castelnuovo's (see \cite[Teorema 1.1]{cc}). Hence $\Delta$ is projectively normal. 
Then also $\Sigma$ is projectively normal, since we have exact sequences:
$$ H^ 1(\pp^ {r+1},\Ical_{\Sigma,\pp^{ r+1}}(n-1))\to H^ 1(\pp^ {r+1},\Ical_{\Sigma,\pp^ {r+1}}(n))
\to H^ 1(\pp^ {r},\Ical_{\Delta,\pp^ r}(n))=0$$
for all $n\geq 1$. 

Now by Lemma \ref{Ciro}  the set $\Gamma$ is contained in  a curve $Z$ of degree $\delta\leq r+1$ which is the intersection of all quadrics containing $\Gamma$. Since $\delta\leq r+1$ the arithmetic genus of $Z$ is $\leq 2$, by Castelnuovo's theorem (cf. \cite[Theorem 3.7]{h}). The same argument  as in the proof of Proposition
\ref {prop:vol1} shows that the intersection of all quadrics containing  $\Sigma$ is a threefold $X$ whose general curve section is $Z$.

Now we want to show that $\delta=r+1$. Suppose otherwise.
If $\delta=r-1$, then $Z$ is a rational normal curve and $h_\Gamma(2)=h_Z(2)=2r-1$, contradicting
$h_\Gamma(2)=2r+1$. So $\delta\geq r$. If $Z$ is rational, then $X$ is not linearly normal (see \cite {eh}) and  $\Sigma$ is not linearly normal either, a contradiction. If $\delta=r$, then $Z$ is a linearly normal smooth elliptic curve, and $h_\Gamma(2)=h_Z(2)=2r$, which again contradicts $h_\Gamma(2)=2r+1$. \end{proof}

In conclusion we have:

\begin{thm}\label{g3}

\begin{enumerate}
\item     If $p_g\ge 8$, then  the intersection of all the quadrics containing the canonical image  $\Sigma$ of $V$ is a linearly normal threefold $X$ of degree $\delta\leq p_g-1$  whose   curve sections have  arithmetic genus $\leq 2$;
\item         if $p_g\geq 14$, then $X$ is a rational normal scroll and the moving part $|D|$  of the pull back $|F|$ on $V$ of the ruling of $X$ is  a base point free pencil of curves of genus 3.
\item         if $p_g\geq 17$, then the  pull back $|F|$ on $V$ of the ruling of $X$ is  a base point free pencil of curves of genus 3.
\end{enumerate}
\end{thm}
\begin{proof} Assertion (i) follows immediately from  Propositions \ref{prop:vol1} and  \ref{prop:vol2}.

 By  Corollary \ref{threefold} and assertion (i), the hypothesis $p_g\geq 14$ ensures that the threefold $X$ is a scroll. Since $V$ is regular by Lemma \ref{q}, $X$ must be a rational normal scroll. We denote by $|F|$ the pull back on $V$ of the ruling of $X$ and we write $|F|=Z+|D|$, where $Z$ is the fixed part and $|D|$ is the moving part of $|F|$.
  The pencil $|D|$  has genus $\ge 3$ because the canonical map of $V$ is birational.
 
 One has  $h_{\Sigma}(4)\leq P_4(V)\leq 19p_g-29$. By using relations of the type \eqref{eq:fact1}, \eqref{fact2} etc. for $X$, its surface and curve sections (cf. \cite{h}, Lemma (3.1) and Remark (3.1.1)), one finds $h_X(4)\ge 20p_g-45$. Hence for  $p_g\ge 17$, one has $h_{\Si}(4)< h_X(4)$ and so there is a quartic hypersurface containing $\Si$ but not $X$. It follows that $F$ has genus 3, $Z=0$ and the pencil $|F|$ is free.
 
 Assume now $14\leq p_g\leq 16$.  
 Arguing as above, one shows  $h_{\Si}(5)\le 31p_g-49$ and $h_X(5)\geq 35p_g-84$. Thus  there is a quintic hypersurface containing $\Si$ and not containing $X$. It follows that $K_VD\le 5$ and therefore  $D^2=0$ by the index theorem. Now $K_VD$ is even by the adjunction formula, hence $K_VD=4$ and $g(D)=3$. 
\end{proof}

\section{Canonical double planes with a free $\Z_2^3$-action}\label{nonbir}
\subsection{The set-up}\label{setup}

In this section we study the following situation:
\begin{itemize}
\item $S$ is a minimal surface of general type;
\item $\pi\colon V \to S$ is an \'etale Galois  cover with Galois group $G$ isomorphic to $\Z_2^3$;
\item the canonical map of $V$ is $2$-to-$1$ onto a rational surface.
\end{itemize}
\begin{rem}{\rm By \cite[Corollary 5.8]{beauville} (cf. \cite[Proposition 4.1]{3chi}) a group $G$ that  acts freely on a surface whose canonical map is of degree 2 on a rational surface is isomorphic to $\Z_2^k$ for some $k$.}
\end{rem}

 Denote by $\iota$ the involution associated with the canonical map of $V$. Then $\iota$ is in the center of $\Aut(V)$, and it induces an involution $\bar{\iota}$ of $S$. Let $\alpha\colon V\to W:= V/\iota$ and $\bar{\alpha}\colon S\to \bar{W}:=S/\bar{\iota}$ be the quotient maps. Notice that $G$ acts also on $W$ and we have  the following commutative  diagram:
\begin{equation} \label{diagram}
\begin{CD}
V @>{\alpha}>>W \\
@V{\pi}VV @VV{\bar{\pi}}V\\
S @>{\bar{\alpha}}>> \bar{W} 
\end{CD}
\end{equation}
where the vertical maps are $G$-covers. 

Let $\tilde{V}$ be the blow up of $V$ at the isolated fixed points of $\iota$. Since the set of isolated fixed points of $\iota$ is $G$-invariant, the action of $G$ on $V$ lifts to an action on $\tilde{V}$, which is again free.  The involution $\iota$ lifts to an involution of $\tilde{V}$ that we denote by the same letter. The quotient  surface $\tilde{W}:=\tilde{V}/\iota$ is smooth and the natural map $\tilde{W}\to W$ is the minimal resolution of the singularities of $W$.

\subsection{The minimal model of $(\tilde{W},G)$.}
Clearly, the group $G$ acts also on $\tilde{W}$. Then (cf.~\cite{cremona1}, \cite{zhang}) $\tilde{W}$ dominates a {\em minimal pair} $(X,G)$, namely there is a smooth rational surface $X$ with the following properties:
\begin{itemize}
\item the group $G$ acts on $X$ and there is a birational $G$-equivariant morphism $\tilde{W}\to X$;
\item every birational $G$-equivariant morphism $X\to X'$, where $X'$ is a smooth rational surface with a $G$-action, is an isomorphism.
\end{itemize}
By \cite[Theorem 4]{zhang} there are the following possibilities for $X$: 
\begin{itemize}
\item[(a)] The group $\Pic(X)^G$ has rank $\ge 2$. In this case, there is a $G$-invariant Mori fibration  $f\colon X \to \pp^1$. In particular, the smooth fibres of $f$ are isomorphic to $\pp^1$ and the singular fibres of $f$ are the union of two irreducible $-1$-curves meeting in a point. 
\item[(b)] The group $\Pic(X)^G$ has rank $1$. In this case, $X$ 
is a Del Pezzo surface.
\end{itemize}
The following elementary observation  is the key of our analysis of the $G$-action on $\tilde{W}$.
\begin{lem}\label{GP}
The stabilizer $G_P<G$ of any point $P\in \tilde{W}$ is cyclic.
\end{lem}
\begin{proof}
Let $P'\in \tilde{V}$ be a point that  maps to $P$. Assume by contradiction that there exist  nonzero elements $g_1\ne g_2\in G_P$. Since $G$ acts freely on $\tilde{V}$, we have  $g_1(P')=\iota (P')=g_2(P')$. It follows that $P'$ is a fixed point of $g_1g_2$ on $\tilde{V}$, a contradiction.
\end{proof}

The statement of  Lemma \ref{GP} holds  also for the minimal model $(X, G)$:

\begin{lem} \label{GPmin}
For every $P\in X$ the stabilizer $G_P<G$ is cyclic.
\end{lem}
\begin{proof} Let $Y$ be  a smooth surface with a faithful $G$-action, and a point $P\in Y$ with nontrivial stabilizer $G_P$. By Cartan's Lemma (see \cite{cartan})  there exist local analytic coordinates in a neighbourhood of $P$ with respect to which the action of $G_P$ is linear. Hence $G_P$ is isomorphic either to $\Z_2$ or to $\Z_2^2$. Denote by $B_2(Y)$ the subset of points of $Y$ with stabilizer isomorphic to $\Z_2^2$.

 By Lemma \ref{GP}, it is enough to show that if $B_2(Y)$ is nonempty, then $B_2(Y')$ is also nonempty for  any pair $(Y', G)$ that dominates $(Y,G)$. We will prove this by induction on the number $k$ of blow-ups of which the map $\epsi\colon Y'\to Y$ is composed.
 
 Let $P\in Y$ be a fundamental point of $\epsi\inv$ and consider the blow up $Y''\to Y$ of all the points in the $G$-orbit of $P$. Then the action of $G$ on $Y$ induces an action on $Y''$ and $\epsi$ is the composition of the $G$-equivariant birational maps $Y'\to Y''$ and $Y''\to Y$. If $P\notin B_2(Y)$, then $B_2(Y'')$ contains the inverse image of $B_2(Y)$, hence it  is nonempty. The claim now follows by induction,  since the number of blowups in $Y'\to Y''$ is strictly smaller than $k$.
 So assume that $P\in B_2(Y)$. Then by Cartan's Lemma, there are generators $g_1$, $g_2$ of $G_P$ and local analytic coordinates $(x,y)$ near $P$ such that $G_P$ acts by:
 \begin{equation}\label{action}
 (x,y)\overset{g_1}{\mapsto}(x,-y);\quad (x,y)\overset{g_2}{\mapsto}(-x,y).
 \end{equation}
 In particular, the fixed loci of $g_1$ and $g_2$ are divisorial near $P$. 
Let $E_P\subset Y''$ be the exceptional divisor over $P$. Then (\ref{action}) shows that $g_1g_2$ acts trivially on $E_P$. Hence the points of intersection of $E_P$ with the strict transforms of the divisorial part of the fixed loci of $g_1$ and $g_2$ are in $B_2(Y'')$. The claim now follows by induction as in the previous case.
\end{proof}
\begin{prop} \label{minimal}\begin{enumerate}
\item The surface  $X$ is a $G$-invariant conic bundle;
\item the subgroup $H<G$ of the elements that map every fibre of $f$ to iself is cyclic generated by an element $h\ne 0$;
\item  the fixed locus of $h$ is a smooth irreducible bisection $D_0$  of $f$ of genus $g(D_0)\equiv 1\mod 4$. There are $2g(D_0)+2$ singular fibres of $f$.  The divisor $D_0$ intersects every singular fibre of $f$ at the singular point and $h$ exchanges the two components of each singular fibre;
\item if $g\in G\setminus H$, then the fixed locus of $g$ consists  of two pairs of points, lying on two smooth fibres of $f$.

\end{enumerate}
\end{prop}
\begin{rem}{\rm An involution $h$ as in Proposition \ref{GPmin} is called a {\it De Jonqui\`eres involution} of degree $g(D_0)+2$ (cf. \cite[\S 2]{cremona1}\ ). Property (iii) implies that $X$ is a minimal model  for $h$.}
\end{rem}

\begin{proof} 
We divide the proof into several steps. 

\medskip

\underline{Step 1:} {\em If $X$ is a Del Pezzo surface and   $\Pic(X)^G$ has  rank 1, then  $K^2_X=8$  does not occur.} \newline
Assume by contradiction that this is the case. 
Then  $X$ is either $\pp^1\times \pp^1$ or $X$ is the blow up of $\pp^2$ in a point. If $X$ is equal to $\pp^1\times\pp^1$ then the condition that $\Pic(X)^G$ has rank 1 implies that $G$ contains at least one (and hence  four)  elements  that exchange the two copies of $\pp^1$. An element $g$ of this type has the form  $(x, y)\mapsto(a\inv y, ax)$, where $a\in \Aut(\pp^1)$.   The  fixed locus of $g$   is the graph of the automorphism $a$ of $\pp^1$, and therefore it is a curve of type $(1,1)$.  Since two such curves have nonempty intersection, we have a contradiction to Lemma \ref{GPmin}. If $X$ is the blow up of the plane in a point, then every automorphism of $X$ maps the unique $-1$-curve of $X$ to itself. Since the canonical class is also preserved by every automorphism,  $G$ acts trivially on $\Pic(X)$, which has rank 2, contradicting the assumptions.
\medskip

\underline{Step 2:} {\em If  $X$ is a Del Pezzo surface and   $\Pic(X)^G$ has  rank 1, then there is precisely one element $h\in G\setminus\{0\}$ such that the fixed locus of $h$ contains a divisor $D_h$.}\newline
 Assume first that the fixed locus of every nonzero element of $G$ has dimension 0. The quotient surface $Y:=X/G$ has only nodes as singularities and we have $K^2_X=8K^2_Y$, hence we have $K^2_X=8$,  contradicting Step 1. So there is at least an element $h$ such that the divisorial part $D_h$ of the fixed locus of $h$ is nonempty. Since $D_h$ is invariant under the group $G$, there exists a positive rational number $r$ such that $D_h$ is numerically equivalent to $-rK_X$. Hence $D_h$ is ample. Since this argument applies to the divisorial part of the fixed locus of any nonzero $g\in G$, by Lemma \ref{GPmin} the fixed locus of every nonzero $g\ne h$ is a finite set.
 \medskip

\underline{Step 3:} {\em The case $X$ Del Pezzo and $\Pic(X)^G$ of rank 1 does not occur, hence there is a $G$-invariant Mori fibration $f\colon X\to\pp^1$.}
 
By Step 2, there is a subgroup $G_1<G$ of order 4 such that the fixed locus of every element of $G_1$ is a finite set. (Recall that by (\ref{TFrat}) the fixed locus of a nonzero element of $G_1$ consists of 4 points). By Cartan's Lemma the quotient surface $Y:=X/G_1$ has only nodes as singularities,  and we have $K^2_X=4K^2_Y$, hence  $K^2_X$ is divisible by 4. Then we have $K^2_X=4$ by Step 1. Notice that in this case $K_X$ is not divisible in  $\Pic(X)$, since $X$ is not minimal. 
It follows that the curve $D_h$ is linearly equivalent to $-rK_X$, where $r>0$ is an integer. Let $\nu$ be the number of isolated fixed points of $h$. Formula \ref{TFrat} gives: $$4=\nu-K_XD_h=\nu+rK_X^2=\nu+4r,$$ i.e. $\nu=0$, $r=1$. Hence the quotient surface $T:=X/h$ is  smooth.  If we denote by $p\colon X\to T$ the quotient map, then the adjunction formula gives $p^*K_T=2K_X$, hence $T$ is a Del Pezzo surface with $K^2_T=8$. The Picard number of $T$ is 2, hence the group $\Pic(X)^h$ has rank 2. 
 The group $G/\!\!<h>\!$ acts on $T$. Let $\ga\in G/\!\!<h>\!$ be a nonzero element and denote by $t_{\ga}$ the trace of $\ga$ in $H^2(T,\C)$. Then, since the fixed locus of $\ga$ consists of 4 points,  the topological trace formula (\ref{TTF}) gives $2+t_{\ga}=4$, namely $\ga$ acts on $H^2(T,\C)$ as the identity.  It follows that $\Pic(X)^G$ has rank 2, contradicting the assumptions. 
\smallskip 

 Hence we have proven that $\Pic(X)^G$ has rank greater than 1, and therefore we conclude that $X$ has the structure of a $G$-invariant conic bundle, proving (i).
 \medskip
 
\underline{Step 4:} {\em  $X=\pp^1\times \pp^1$ does not occur.}\newline
 Assume that this is the case.  By Step 3, the group $G$ acts trivially on $\Pic(X)$, hence in particular it preserves the projections $p_i\colon\pp^1\times\pp^1\to \pp^1$, $i=1,2$. Thus  every element of $G$ is of the form $(x,y)\mapsto(ax,by)$, where $a, b\in \Aut(\pp^1)$ have order 2. Consider the homomorphism $\psi_1\colon G\to \Aut(\pp^1)$ induced by the projection $p_1$. Since $\Aut(\pp^1)$ does not contain a subgroup isomorphic to $\Z_2^3$, the kernel of $\psi_1$  contains a nontrivial  element $g_1$  of the form $(x,y)\mapsto (x, by)$. By the same argument, $G$ contains also a nontrivial element 
 $g_2$ of the form $(x,y)\mapsto (ax, y)$. Then the fixed loci of $g_1$ and $g_2$ intersect, contradicting Lemma \ref{GPmin}.
 \medskip
 
\underline{Step 5:} {\em  The subgroup  $H< G$ consisting of the elements that map the general fibre of $f$ to itself is isomorphic to $\Z_2$.  The quotient group $G/H$ acts freely on the set of singular fibres of $f$.}\newline
By definition, the group $H$ acts faithfully  on the general fibre of $f$, which is a $\pp^1$,  and $G/H$ acts faithfully on the base of the fibration $f$, which is also a $\pp^1$.  
 Since $\pp^1$ has no faithful $\Z_2^3$-action,  both $H$  and $G/H$ have  order at most 4. In particular, $H$ has either  order 2 or 4. Let $h \in H$ be a nonzero element. The fixed locus of $h$ on a general fibre of $f$ consists of two points. Hence the fixed locus of $h$ contains a bisection  $D_0$ of the fibration $f$. If $H$ has order 4, then $f$ has no reducible fibre since otherwise the singular point of the fibre would be fixed by every element of $H$, contradicting Lemma \ref{GPmin}. Hence $f\colon X\to\pp^1$ is a $\pp^1$-bundle. Since by Lemma \ref{GPmin} the fixed loci of the three nonzero elements of $H$ are disjoint, $X$ contains three disjoint bisections. This is possible only if $X=\pp^1\times \pp^1$, contradicting Step 4. So $H$ has order 2.
 
 Assume by contradiction that the action of $G/H$ on the set of singular fibres of $f$ is not free. Then there is a singular fibre of $f$ that is mapped to itself by a subgroup $H_0\subset G$ that contains properly $H$. It follows that the singular point of the fibre  is fixed by $H_0$, contradicting Lemma \ref{GPmin}.\medskip
 
\underline{Step 6:} {\em Let $h$ be the generator of $H$. Then the fixed locus of $h$ is a smooth bisection $D_0$ of $f$,  which intersects the reducible fibres of $f$ at the singular point and the smooth fibres at two distinct points. The two irreducible components of a singular fibre are exchanged by $h$.}\newline
As we remarked in Step 5,  the fixed locus of $h$ contains a bisection  $D_0$ of the fibration $f$.
The divisorial part $D_h$ of the fixed locus of $h$ is smooth, hence $D_h$ does not contain any fibre of $f$. We claim that if $A_1+A_2$ is a reducible fibre of $f$, then $h$ exchanges $A_1$ and $A_2$ and, as a consequence, $D_h$ intersects $A_1+A_2$ at the common point  of $A_1$ and $A_2$. 
Assume by contradiction that $h(A_1)=A_1$.  By Step 5 the group $G/H$ acts freely on the set of reducible fibres of $f$, hence the  $G$-orbit of $A_1$ is a disjoint union of four $-1$-curves,  contradicting the minimality of $(X,G)$. 

Now, in order to show that the fixed locus of $h$ coincides with $D_0$,  it is enough to prove that $D_0$ meets every smooth fibre of $f$ at two distinct points. Assume by contradiction that there exists a point $P\in X$ such that $D_0$ is tangent at $P$ to the fibre $C$ of $f$ containing $P$. Then $h$ acts as the identity on the space $T_PD_0=T_PC$. By Cartan's Lemma this contradicts the fact that $h$ maps $C$ to itself but it does not fix $C$ pointwise.
\medskip

\underline{Step 7:} {\em The curve $D_0$ is irreducible and $g(D_0)-1$ is divisible by 4. The fibration  $f$ has precisely  $2g(D_0)+2$ reducible fibres. }\newline
By Lemma \ref{GPmin}, the quotient group $G/H$ acts freely on $D_0$. 
Assume that $D_0$ is reducible.  In this case $f$ is a $\pp^1$-bundle by Step 6, and $D_0$ is the disjoint union of two sections of $f$.  Since $G/H$ has order 4 there exists a nonzero $\ga\in G/H$ that maps each component of $D_0$ to itself. Since these components are rational curves, this contradicts the fact that $G/H$ acts freely on $D_0$. Hence $D_0$ is irreducible and  $g(D_0)-1$ is divisible by 4 by the Hurwitz formula. By Step 6 the singular fibres of $f$ are in $1$-to-$1$ correspondence with the ramification points of the double cover $D_0\to\pp^1$ induced by $f$, hence there are $2g(D_0)+2$ of them by the Hurwitz formula. \medskip

\underline{Step 8:} {\em If $g\notin H$, then the fixed locus of $g$ consists of 4 points.}
Consider now an element $g\in G\setminus H$. No reducible fibre of $f$ can be mapped to itself by $g$, since otherwise the singular point of the fibre would be fixed by both $g$ and $h$, contradicting Lemma \ref{GPmin}. Hence the fixed locus of $g$ is contained in two smooth fibres of $f$ and,  again by Lemma \ref{GPmin}, it consists of two   pairs of points, one  on each  fibre. 
\end{proof}
\begin{rem} {\em The minimal pair $(X,G)$ can be constructed as follows. Let $\psi\colon Z\to\pp^1$ be a geometrically ruled rational surface and let $C\subset Z$ be a smooth bisection of genus $g(C)>0$ such that $C$ is divisible by 2 in $\Pic(Z)$. Let $P_1, P_2, P_3$ be points of $C$ such that  for $i=1,2,3$ the curve $C$ is tangent  to the fibre of $\psi$ in $P_i$ (there are $2g(C)+2\ge 4$ such points by the Hurwitz formula). Let $\epsi\colon Z'\to Z$ be the surface obtained by blowing up first  the points $P_1,P_2,P_3$ and then for $i=1,2,3$ the intersection point of the exceptional curve over $P_i$ with the strict transform of the fibre of $\psi$ through $P_i$. Let $\psi'\colon Z'\to \pp^1$ be the fibration induced by $\psi$. The map $\psi'$ has  3 reducible fibres $F_1, F_2, F_3$, which are the pull back of the fibres of $\psi$ containing $P_1,P_2,P_3$, respectively. For every $i=1,2,3$ we write $F_i=2E_i+N_{2i-1}+N_{2i}$, where $E_i$ is a $-1$-curve, $N_{2i-1}$ is the strict transform of the fibre of $\psi$ through $P_i$ and $N_{2i}$ is the strict transform of the exceptional curve over $P_i$. The curves $N_1, \dots N_6$ are disjoint  $-2$-curves. One has:
$$\epsi^*C=C'+\sum_{i=1}^{3}(2E_i+N_{2i}),$$
where $C'$ is the strict transform.

Let $g_1, g_2,h$ be generators of $G$ and let $\chi_1,\chi_2,\chi_3\in G^*$ be the dual basis of characters.
Write $g_3=g_1g_2$ and set:

\begin{gather*} D_{g_1}=N_1,\  D_{g_2}=N_3,\  D_{g_3}=N_5,\  D_h=C'\\
D_{g_1h_1}=N_2,\  D_{g_2h_2}=N_4, \ D_{g_3h_3}=N_6.
\end{gather*}
The reduced fundamental relations (cf. \cite[Proposition 2.1]{ritaabel}) for a $G$-cover with branch divisors $D_g$ as above are the following:
$$2L_{\chi_1}\equiv D_{g_1}+D_{g_1h}+D_{g_3}+D_{g_3h}=N_1+N_2+N_5+N_6,$$
$$2L_{\chi_2}\equiv D_{g_2}+D_{g_2h}+D_{g_3}+D_{g_3h}=N_3+N_4+N_5+N_6,$$
$$2L_{\chi_3}\equiv D_{h}+D_{g_1h}+D_{g_2h}+D_{g_3h}=C' +N_2+N_4+N_6.$$
It is not difficult to check that the above relations admit a solution, which is unique since $Z'$ is regular. Hence by \cite[Proposition 2.1]{ritaabel} there exists a $G$-cover $X'\to Z'$ with reduced building data $L_{\chi_i}$, $D_g$ as above. The surface $X'$ is smooth by \cite[Proposition 3.1]{ritaabel} and it  is rational, since the pull back of a general fibre of $\psi'$ is the union of $4$ disjoint smooth rational curves. For $i=1, \dots 6$ the inverse image  of $N_i$ is the disjoint union of four $-1$-curves. Let $X$ be the surface obtained from $X'$ by contracting all these $-1$-curves. Then the $G$-action on $X'$ descends to $X$ and $(X,G)$ is a minimal model as in the statement of Proposition \ref{minimal}}.\end{rem}
\subsection{Numerical invariants}
Consider now diagram (\ref{diagram}).
The surfaces $W$ and $\bar{W}$ are rational surfaces. The singularities of $W$ and $\bar{W}$ are nodes, which are the images of the isolated fixed points of $\iota$, respectively $\bar{\iota}$.
We let $\rho$, respectively $\bar{\rho}$,  denote the number of isolated fixed points of $\iota$,  respectively $\bar{\iota}$ and we let $B$, respectively $\bar{B}$, denote the branch divisor of $\alpha$, respectively $\bar{\alpha}$. Recall that $B$ and $\bar{B}$ are smooth curves.

Set $G_0:=G\times\!\!<\iota>\,\subseteq\Aut(V)$. For every $\gamma\in G_0$, $\ga\ne\iota$, denote by  $n_{\ga}$ the number of isolated fixed points of $\ga$ on $V$ and denote by $D_{\ga}$ the divisorial part of the fixed locus of $\ga$ on $V$.
Notice that for $g\in G$ we have $n_{g}=0$, $D_g=0$. 
Similarly, we denote by $\De_{g}$ the divisorial part of the fixed locus of $g$ on $W$ and by $\nu_{g}$ the number of isolated fixed points of $g$ on $W$.
By the commutativity of diagram \ref{diagram}, we have:
\begin{equation}
 \bar{B}=\bar{\pi}(B)+ \sum_{g\ne 1}\bar{\pi}(\Delta_g)
 \end{equation}
Analogously, we have the following formula for the number of isolated fixed points:
\begin{equation}\label{rho1}
\bar{\rho}=\frac{1}{8}\rho+\frac{1}{8}\sum_{g\ne 1} n_{g\iota}
\end{equation}

Recall (see \S \ref{setup}) that $W$ has canonical singularities. In particular, it is Gorenstein and $2$-factorial, hence  the intersection number $K_WD$ is defined for any divisor $D$ of $W$.

Using the analysis of the previous subsection we can prove:
\begin{prop}\label{Dg}
Let $1\ne g\in G$. Then: 
\begin{enumerate} 
\item if $C\subset W$ is a component of $\De_g$, then $K_WC\ge 0$, and equality holds if and only if $C$ is a $-2$-curve;
\item $n_{g\iota}=2\nu_g=8+2K_W\De_g$;
\item for $g=h$ {\rm (cf. Proposition \ref{minimal})}, the divisor  $\De_h$ contains an hyperelliptic curve $D$ with  $g(D)\equiv 1 \mod 4$.
\end{enumerate} 
\end{prop}
\begin{proof} 
Recall that, since the $G$-action on $W$ is induced by a free action on $V$, $G$ acts freely on the branch locus of $\iota$, and in particular on the nodes of $W$. Hence, if we denote by $\De$ the set of points $P$ of $W$ such that the stabilizer $G_P<G$ is non trivial, then  $\alpha$ is \'etale above $\Delta$. In particular, for every curve $C\subset \De$ one has $2K_WC=K_V\alpha^*C$, and statement (i) follows since $V$ is minimal of general type.
 
The equality $n_{g\iota}=2\nu_g$ is a consequence of the previous remark and of $n_g=0$. In addition, we have $\nu_g=4+K_W\De_g$ by (\ref{TFrat}) and Remark \ref{remsing}. 

Finally statement (iii) follows directly by Proposition \ref{minimal}.
\end{proof}
\begin{cor}\label{8} 
One has $\bar{\rho}\ge 8$;
\end{cor}
\begin{proof} The result follows  from formula (\ref{rho1}) and Proposition  \ref{Dg}.
\end{proof}
\bigskip

Let $\eta\colon S'\to S$ be the blow up of the isolated fixed points of
 $\bar{\iota}$ and let $Y\to \bar{W}$ be the minimal desingularization.
Note that  $\bar B$, which is a curve on $\bar W$, can be seen as a curve
 on $Y$.  Let
$p\colon S'\to Y$ be the map induced by $S\to \bar{W}$. The map $p$ is a
 flat double cover of smooth surfaces, given by an equivalence relation
 $2L\equiv \bar{B}+N$, where $N$ is the
 exceptional divisor of $Y\to\bar{W}$.

By standard formulae for double
 covers  one has:
\[
 H^0(K_{S'})=H^0(K_Y)\oplus H^0(K_Y+L),
 \]
 \[
 H^0(2K_{S'})=H^0(2(K_Y+L))\oplus H^0(2K_Y+L).
 \]
In the above decompositions  the first summand corresponds to ``even'' sections and the second one to
 ``odd'' sections.
 
 The dimension of $H^0(2K_Y+L)$ can be computed as follows (cf. \cite[\S 3]{ccm}  for  similar computations):
\begin{prop}\label{2K+L} One has:
\begin{enumerate}
\item $K^2_S\equiv \bar{\rho}  \mod 2$;
\item $2h^0(2K_Y+L)=6+K^2_S-\bar{\rho}-2\chi(S)$.
\end{enumerate}
\end{prop}
\begin{proof}
Since $K_{S'}=p^*(K_Y+L)$, we have:
\begin{equation}\label{K2}
K^2_S-\bar{\rho}=K^2_{S'}=2(K_Y+L)^2,
\end{equation}
proving statement (i).

Computing $\chi(S)=\chi(S')$ by means of  the standard formulae for double covers, we get:
\begin{equation}\label{chi}
L(K_Y+L)=2\chi(S)-4.
\end{equation}
By the Kawamata-Viehweg vanishing theorem, we have $h^0(2K_Y+L)=\chi(2K_Y+L)$, hence
\begin{equation}
h^0(2K_Y+L)= 1+\frac{1}{2}(2K_Y+L)(K_Y+L)=1+(K+L)^2-\frac{1}{2}L(K_Y+L)
\end{equation}
By (\ref{K2}) and (\ref{chi}), this is equivalent to statement (ii).
\end{proof}

\section{The fundamental group of surfaces with $K^2=3\chi-1$}\label{sec:fund}
The aim of this section is to prove the following:
\begin{thm}\label{main1}
Let $S$ be a minimal complex surface of general type such that $K^2_S=3\chi(S)-1$.
If the group $\pionealg(S)$ has order 8, then $\chi(S)=1$. 
\end{thm}

\medskip

We fix the notation and assumptions that will hold throughout all the section.

We consider a smooth minimal complex surface of general type $S$ such that:
\begin{itemize}
\item $K^2_S=3\chi(S)-1$;
 \item $G:=\pionealg(S)$ has order 8.
\end{itemize}

We denote by $\pi\colon V\to S$  the (algebraic) universal cover of $S$. The surface $V$ has invariants $K^2_V=8K_S^2$, $\chi(V)=8\chi(S)$, $q(V)=0$, $p_g(V)=\chi(V)-1$. In particular one has $K^2_V=3\chi(V)-8=3p_g(V)-5$. We start by analyzing the canonical map of $V$. 

\begin{prop}\label{canV} Set $r:=\chi(V)-2$,  let $\Si$ be the canonical image of $V$ and let  $\fie\colon V\to\Si\subset \pp^r$ be the canonical map. Then there are the following possibilities:
\begin{enumerate}
\item $\fie$ is a degree 2 map and $\Si$ is a rational surface;
\item $\fie$ is birational and the system $|K_V|$ is free.
\end{enumerate}
\end{prop}
\begin{proof}
We start by showing that $\fie$ is not composed with a pencil. This also follows from \cite{ho}  but here we give a direct proof. Assume by contradiction that this is the case and write $|K_V|=Z+|rF|$, where $|F|$ is a pencil. Then:
\[3r-2=K^2_V\ge rK_VF \]
which implies $K_VF\le 2$. Now the index theorem gives $F^2=0$, hence  by the adjunction formula one has 
$K_VF=2$ and $|F|$ is a $G$-invariant genus 2 pencil, contradicting Lemma \ref{free}.
So $\Si$ is a surface. 

Let $d$ be the degree of $\fie$. Write $|K_V|=Z+|M|$, where $Z$ is the fixed part and $|M|$ is the moving part.
If $Z\ne 0$, then $MZ\ge 2$ by the $2$-connectedness of canonical divisors. 
We have the following chain of inequalities:
\begin{gather}\label{chain}
3r-2=K^2_V=K_VM+K_VZ\ge K_VM=\\=M^2+MZ\ge M^2\ge d\deg\Si.\nonumber
\end{gather}
Since a nondegenerate surface of $\pp^r$ has degree $\ge r-1$, (\ref{chain}) gives  $d\le 3$.

Assume   $d=3$.  In this case  (\ref{chain})  gives  $\deg\Si=r-1$ and $MZ\le 1$. So  we have $MZ=0$ and $Z=0$. The relation $K^2_V=3\deg\Si+1$ implies that $|K_V|$ has precisely one base point, which must be fixed by every element of $G$, contradicting the fact that $G$ acts freely on $V$. So $d=3$ does not occur.

Assume $d=2$. Then (\ref{chain}) gives  $\deg\Si\le \frac{3}{2}r-1<2r-2$  (recall that $r\ge 6$) and by \cite[Lemme 1.4]{beauville} the surface $\Si$ is ruled. In addition, $\Si$ is regular, since it is dominated by $V$,  and therefore it is a rational surface.

Assume $d=1$. Then by the Castelnuovo inequality (\cite[Remarques 5.6, 1)]{beauville}) we have $\deg\Si\ge 3r-4$. Hence (\ref{chain}) gives $MZ\le 2$. In addition, if $MZ=2$ then $K_VZ=0$ and $Z^2=-2$, hence $Z$ is a connected configuration of $-2$-curves  and the group $G$ maps $Z$ to itself.   Considering the action of  $G$ on the Dynkin diagram associated  to $Z$, it is easy to check that every element of $G$ either maps a component of $Z$ to itself or it exchanges two components of $Z$ that intersect in a point. In either case, this contradicts the assumption that  $G$ act freely on $V$. This proves $MZ=0$, hence $Z=0$ and $|K_V|$ has no fixed part.

Assume that the base locus ${\mathcal  B}$ of $|K_V|$ is not empty.  Since $G$ maps  ${\mathcal  B}$ to itself,  the cardinality of ${\mathcal  B}$ is equal to $8m$ for some positive integer $m$. 
Then we have $K^2_V-8m=3r-2-8m\ge \deg\Si$,  contradicting $\deg\Si\ge 3r-4$.
\end{proof}

The two possibilities  for the canonical map of $V$ given in Proposition \ref{canV} correspond to the situations studied in sections \ref{sec:cansurf}} and \ref{nonbir}. First of all we rule out case (i), namely we prove:
\begin{prop}\label{birV}
The canonical map of $V$ is birational.
\end{prop}
\begin{proof}
Assume by contradiction that the canonical map of $V$ is not birational. Then by Proposition \ref{canV} the canonical map of $V$ is of degree 2 onto a rational surface. By \cite[Cor. 5.8]{beauville} (cf. \cite[Prop. 4.1]{3chi}) the group $G$ is isomorphic to $\Z_2^3$  and we can apply the results of section \ref{nonbir}. Throughout the proof we use freely all the notation of  section \ref{nonbir}.

Recall that  we have $\bar{\rho}\ge 8$ by Corollary \ref{8}. On the other hand, by Proposition  \ref{2K+L}, (ii),  we have $\bar{\rho}\le K^2_S-2\chi(S)+6=\chi(S)+5$. Hence $\chi(S)\ge 3$, and $p_g(S)\ge 2$.

Let $D$ be the smooth curve of $W$ defined in Proposition \ref{Dg}, (iii).
The intersection number $K_WD$ is strictly positive by Proposition \ref{Dg}. In addition, as we observed in the proof of Proposition \ref{Dg}, one has $2K_WD=K_V\alpha^*D$. Since the divisor $\alpha^*D$ is clearly $G$-invariant and $G$ acts freely on $V$, $K_V\alpha^*D$ is divisible by 8 and we may write 
$K_WD=4\epsilon$, for an integer  $\epsilon>0$. Set: $$4\ga:=\sum_{g\in G\setminus<h>}K_W\De_g +K_W(\De_h-D)+\frac{1}{2}\rho.$$
 
Note that since $D\le \De_h $, $\De_h-D$ is effective. 
  By Proposition \ref{Dg} (i), $K_W$ is nef on the components of $\De_g$ and thus $\ga\ge 0$.

By formulae (\ref{rho1}) and Proposition \ref{Dg},  one has: \begin{equation}
\bar{\rho}=\epsilon+7+\gamma
\end{equation}
So the statement of Proposition \ref{2K+L} can be rewritten as follows:
\begin{equation} \label{KA}
2h^0(2K_Y+L)=\chi(S) -2-\epsilon-\gamma=p_g(S)-1- \epsilon-\gamma
\end{equation}
In particular we conclude that $\epsilon\leq p_g(S)-1$.

Let  $\bar{D}:=\bar{\pi}(D)$ and denote again by $\bar{D}$ the corresponding curve in $Y$.  Recalling that $K_Y+L$ pulls back on $S'$ to the canonical bundle and chasing through diagram   (\ref{diagram}),  one checks  that $(K_Y+L)\bar{D}=\epsilon$. Since $g(\bar{D})\geq 1$, we obtain $h^0(\bar{D},\OO_{\bar{D}}(K_Y+L))\leq \epsilon$. 

Consider the restriction map: 
$$H^0(Y,K_Y+L)\to H^0(\bar{D},\OO_{\bar{D}}(K_Y+L)).$$
Taking into account that $h^ 0(K_Y+L)=p_g(S)$,
we see that  we can write $K_Y+L=\bar{D}+A$, where $h^0(Y, A)\geq
p_g(S)-\epsilon>0$.

Since $g(\bar{D})\geq 1$ and the surface $Y$ is regular,  one has $h^0(Y, K_Y+\bar{D})\geq 1$. 
Since  $2K_Y+L=K_Y+\bar{D}+A$, we conclude that $h^0(Y, 2K_Y+L)\geq p_g(S)-\epsilon$, which contradicts \eqref{KA}.
\end{proof}
We notice the following corollary of Proposition \ref{canV} and  Proposition \ref{birV}:
\begin{cor}
Let $S$ be a numerical Campedelli surface, namely a smooth minimal surface of general type with $K^2_S=2$ and $p_g(S)=0$. Assume that $\pionealg(S)$ has order 8 and let $V$ be the (algebraic) universal cover of $S$. Then the system $|K_V|$ is free and the  canonical map of $V$ is birational.
\end{cor}

\medskip Finally we can give the:

\begin{proof}[Proof of Theorem \ref{main1}] The system $|K_V|$ is free and the  canonical map of $V$ is birational by Proposition \ref{canV} and Proposition \ref{birV}. The invariants of $V$ satisfy $K^2_V=3p_g(V)-5$. Assume that  $\chi(S)\ge 2$. In this case $\chi(V)\ge 16$, hence $p_g(V)\ge 15$ and, by Theorem \ref{g3} (ii),  the intersection of all quadrics containing the canonical image $\Si$ of $V$ is a rational normal scroll of dimension 3 and the pull back on $V$ of the ruling of $X$ is a free pencil $|F|$ of curves of genus 3. Since, excepting the cone over the smooth  quadric of $\pp^3$,  a rational normal scroll of has only one ruling,  the pencil $|F|$ is intrinsically attached to $V$.  Therefore $|F|$  is  $G$-invariant, contradicting Lemma \ref{free}. This proves that $\chi(S)=1$.
\end{proof}

\section*{Acknowledgements}
 
The second author is a member of the Center for Mathematical
Analysis, Geometry and Dynamical Systems UTL  and the first and third author are  members of G.N.S.A.G.A.-I.N.d.A.M.  This research was partially supported by the italian  project ``Geometria sulle
variet\`a algebriche'' (PRIN COFIN 2004) and by FCT (Portugal) through program POCI 2010/FEDER and Project
POCTI/MAT/44068/2002.

\end{document}